\documentclass{amsproc}


\newtheorem{theorem}{Theorem}
\newtheorem{lemma}[theorem]{Lemma}
\newtheorem{proposition}[theorem]{Proposition}
\newtheorem{corollary}[theorem]{Corollary}
\newtheorem{question}[theorem]{Question}

\theoremstyle{definition}
\newtheorem{example}[theorem]{Example}

\def\N{{\mathbb N}}

\def\C{{\mathbb C}}

\def\Spec#1{\Phi _ {#1}}
\def\skipaline{\vskip 0.4cm}
\def\skiphalfaline{\vskip 0.2cm}
\def\Proof{{\bf Proof.} }
\def\QED{\hfil\qed\newline\skiphalfaline}
\begin{document}
\title{Countable linear combinations of characters on commutative Banach algebras}
\author{J. F. Feinstein}
\address{School of Mathematical Sciences\\
University of Nottingham\\
University Park\\
Nottingham, NG7~2RD\\
UK.}

\email{Joel.Feinstein@nottingham.ac.uk}
\maketitle

\begin{abstract}
An elegant but elementary result of Wolff from
1921 (\cite{WO}), when interpreted in terms of Banach algebras,
shows that it is possible to find a sequence of distinct characters
$\phi_n$ on the disc algebra and an $\ell_1$ sequence of complex numbers
$\lambda_n$, not all zero, such that
$\sum_{n=1}^\infty \lambda_n \phi_n =0.$
We observe that, even for general commutative, unital Banach algebras,
this is not possible if the closure of the countable set
of characters has no perfect subsets.
\end{abstract}

It is well-known that every finite set of
distinct characters on a commutative Banach algebra is
linearly independent.
In this note we discuss related questions (and their answers)
concerning countable linear combinations of characters.

\skipaline
{\bf Notation.} Throughout this note, $\ell_1$ and $\ell_\infty$ will
denote the usual spaces of sequences of complex numbers.
For a Banach algebra $A$, we denote the character space of $A$ by $\Phi_A$.
\skipaline

We begin with the main question that we wish to discuss.

\begin{question}
Do there exist a commutative, unital Banach algebra $A$, a sequence of
distinct characters $(\phi_n) \subseteq \Spec{A}$ and a sequence
$(\lambda_n) \in \ell_1 \backslash \{0\}$, such that (with absolute
operator norm convergence)
\[
\sum_{n=1}^\infty{\lambda_n \phi_n} = 0?
\]
\end{question}

Suppose for the moment that we have found an example of such a commutative, unital Banach algebra $A$ and sequence of
characters $(\phi_n) \subseteq \Spec{A}$.
Then, for distinct positive integers $m$ and $n$, the set $U_{m n} := \{a \in A: \phi_m(a) \neq \phi_n(a)\}$ is a dense open subset of $A$. By the Baire category theorem, the intersection of all of these sets $U_{m n}$ is dense in $A$. In particular, there exists an $a \in A$ such that the complex numbers $\phi_n(a)$ $(n \in \N)$ are all distinct.
Setting $\alpha_n = \phi_n(a)$ $(n \in \N)$, and considering the elements $a^k \in A$ $(k \in \N)$, we would then obtain a positive answer to the following question about complex sequences.

\begin{question}
Do there exist sequences $(\lambda_n) \in \ell_1\backslash \{0\}$ and
$(\alpha_n) \in \ell_\infty$ with all the $\alpha_n$ distinct and such that
\[
\sum_{n=1}^\infty{\lambda_n{\alpha_n}^k} = 0~~(k=0,1,2,\dots)?\eqno(1)
\]
\end{question}

Conversely, given a positive answer to Question 2, we may obtain a uniform algebra $A$ and a sequence of characters on $A$ satisfying the conditions of Question 1. Indeed, by scaling, we may assume that we have sequences $(\lambda_n)$ and $(\alpha_n)$ satisfying the conditions of Question 2 and such that $|\alpha_n|<1$ for all $n \in \N$.
Take $A$ to be the disc algebra, and let $\phi_n$ be evaluation at the point $\alpha_n$ of the open unit disc. Since the polynomials are dense in the disc algebra, it follows easily that the conditions of Question 1 are satisfied.
\skiphalfaline
Question 2 is discussed in some detail in \cite{BSZ}. They first observe that
(1) is equivalent to each of the following two conditions:
\[
\sum_{n=1}^\infty {\lambda_n \exp ({\alpha_n}z)} = 0~~~~(z \in
\C);\eqno(2)
\]
\[\sum_{n=1}^\infty \frac{\lambda_n}{z - \alpha_n} = 0~~~~ (z \in \C, |z|>\sup_{n}|\alpha_n|).\eqno(3)\]
The answer to Question 2 is yes. Indeed, as noted in \cite{BSZ}, the first
example of such a pair of
sequences $(\lambda_n)$, $(\alpha_n)$ (satisfying (3)) was found
by Wolff in \cite{WO}. His example is so elementary and elegant that we include it
here, along with some extra details.

{\bf Notation.} Let $a\in\C$ and let $r>0$. Then we denote by
$\Delta(a,r)$,
$\bar{\Delta}(a,r)$ respectively the open disc and the closed disc in
$\C$
with centre at $a$ and radius $r$. We denote Lebesgue area measure on
$\C$ by $m$.
\skipaline

\begin{example} \cite{WO} Let $D_1$ be the closed unit disc in $\C$,
$\bar{\Delta}(0,1)$.
We may choose
sequences $(\alpha_n)_{n=2}^\infty \subseteq \Delta(0,1)\backslash\{0\}$
and $r_n > 0$ such that the closed
discs $\bar{\Delta}(\alpha_n,r_n)$ are pairwise disjoint subsets of
$\Delta(0,1)$,
and such that
\[m\left(D_1\backslash\bigcup_{n=2}^\infty
\bar{\Delta}(\alpha_n,r_n)\right)=0\,.\]
Set $\alpha_1=0$ and, for
$n=2,3,\dots$, set $D_n=\bar{\Delta}(\alpha_n,r_n)$. Set
$\lambda_1 = -\pi = -m(D_1)$
and, for $n=2,3,\dots$ set
$\lambda_n=\pi {r_n}^2 = m(D_n)$, so that
\[\sum_{n=1}^\infty{|\lambda_n|} = 2\pi < \infty\,.\]
Then it is clear that, for any complex-valued function $f$ which
is harmonic on a neighbourhood of $D_1$, we have
\[\pi f(\alpha_1) = \sum_{n=2}^\infty{\pi {r_n}^2 f(\alpha_n)}\]
and so
\[\sum_{n=1}^\infty{\lambda_n f(\alpha_n) = 0}. \eqno (4)\]
It is now easy to see that the sequences $(\alpha_n)$, $(\lambda_n)$
have the required properties.
To check (2), for example, let $z\in\C$. Set $f(w) = \exp(z w)$,
so that $f$ is an entire function of $w$.
Applying (4) to $f$ gives
\[\sum_{n=1}^\infty{\lambda_n \exp(\alpha_n z)} = 0\,,\]
as required.
\end{example}
\par
Note that this shows that there is a non-zero annihilating measure for the
disc
algebra which is concentrated on a countable set, given by the countable linear
combination of point masses
\[
\sum_{n=1}^\infty{\lambda_n \delta_{\alpha_n}}.
\]
Also, as discussed earlier, the corresponding countable linear combination of distinct evaluation
characters $\varepsilon_{\alpha_n}$ on the disc algebra is zero: with absolute convergence in
operator norm,
\[\sum_{n=1}^\infty{\lambda_n \varepsilon_{\alpha_n}} = 0\,.\]

These characters are not in the Shilov boundary of the disc algebra:
the F. and M. Riesz theorem (\cite[Chapter II, 7.10]{G})  shows that that would be impossible
for the
disc algebra.
However, it is possible, in general, for all of the characters
to be
in the Shilov boundary. Indeed, the example above may be used to give such
a sequence of characters for the well-known {\it tomato can algebra\/}
(see \cite[Chapter I, Exercise 12]{G} or \cite[2-8]{Browder})
of functions continuous on a solid cylinder, and analytic on one face.
In this case, the Shilov boundary of the algebra is equal to its character space.
This example is not essential. (The term essential was introduced in \cite{Bear}: see also \cite[2-8]{Browder}).
However, by using de Paepe's construction \cite{dePaepe} in place of the tomato can algebra, we may obtain
such an example where the uniform algebra is essential and where the Shilov boundary is again equal to the character space.
\skiphalfaline
If some extra conditions are placed on the sequence $(\alpha_n)$
then no such sequence $(\lambda_n)$ can exist. For example, in \cite{BSZ}, the
following result appears, which is an elementary consequence of the standard
theory of polynomial approximation.
\par
\begin{proposition}
Let $(\alpha_n)$ be a bounded sequence of distinct complex
numbers.
Suppose that the closure in $\C$ of $\{\alpha_1,\alpha_2,\dots\}$ has no
interior, and
does not separate the plane. Then there is no sequence
$(\lambda_n) \in \ell_1 \backslash\{0\}$ satisfying (1).
\end{proposition}
\par
There are some more general results in \cite{BSZ}. For example, the authors give a full
characterisation of those sequences of distinct complex numbers
$(\alpha_n) \subseteq \Delta(0,1)$ for which there is a sequence
$(\lambda_n) \in \ell_1 \backslash\{0\}$ satisfying (1), under the
further assumption that
the sequence $(\alpha_n)$ has no interior limit points in $\Delta(0,1)$.
Such a sequence $(\lambda_n)$ then exists if and only if almost every
point of the
unit circle is a non-tangential limit of a subsequence of $(\alpha_n)$.
They also showed that this is equivalent to the requirement that the
sequence $(\alpha_n)$ is a {\it dominating sequence}, i.e., every bounded
analytic function $f$ on the unit disc satisfies
\[\sup_{z\in\Delta}|f(z)| = \sup_{n}|f(\alpha_n)|\,.\]
\par
Proposition 4 suggests that if some extra conditions are placed on either
the
sequence or the algebra in Question 1, then the answer may change.
For example, Proposition 4 may be used, along with a
Baire category theorem argument, to prove the following result. We shall
not include the elementary proof, since the result also follows from the
more general result (Theorem 9) below.
\par
\begin{theorem}
Let $A$ be a commutative, unital Banach algebra, and let
$(\phi_n)$ be a sequence of distinct characters on $A$ such that
$\{\phi_1, \phi_2,\dots \}$ is a compact subset of $\Spec{A}$.
Suppose that
$(\lambda_n) \in \ell_1$, and that (with absolute operator norm convergence)
\[\sum_{n=1}^\infty{\lambda_n \phi_n} = 0\,.\]
Then all of the $\lambda_n$ must be zero.
\end{theorem}
\par
An alternative approach that yields rather more is to use the
following easy lemma.
\begin{lemma}Suppose that $A$, $(\lambda_n)$ and $(\phi_n)$ satisfy the
conditions of Question 1. Let $E$ be the closure in $\Spec{A}$ of
$\{\phi_1,\phi_2,\dots\}$. Let $B$ be the uniform closure in $C(E)$ of
the restriction to $E$ of the Gelfand transform of $A$. Then $B$ is
a non-trivial uniform algebra on $E$.
\end{lemma}

\par
\Proof
Clearly $B$ is a uniform algebra on $E$, and $B$ is annihilated by the
non-zero, regular Borel measure $\mu$ on $E$ defined by
\[\mu = \sum_{n=1}^\infty{\lambda_n \delta_{\phi_n}}\,.\]
The result follows immediately.
\QED

This is more than enough to resolve Question 1 for one obvious restricted
class of Banach algebras.
\par
\begin{corollary}
Let $A$ be a commutative, unital Banach algebra such that
the Gelfand transform $\hat{A}$ is self-adjoint. Then there are no sequences
$(\lambda_n)$, $(\phi_n)$ satisfying the conditions of Question 1.
\end{corollary}
\Proof
The result follows immediately from Lemma 6 and the Stone-Weierstrass theorem.
\QED
We now return to conditions on the sequence of characters $(\phi_n)$.
To help us we quote the following result of Rudin.
\par
\begin{proposition}
\cite[Theorem 4]{R}. Let $X$ be a compact space which has no
perfect subsets. Then there are no non-trivial uniform algebras
on $X$.
\end{proposition}
\par
We can now give a stronger version of Theorem 5.
\par
\begin{theorem}
 Let $A$ be a commutative, unital Banach algebra, and let
$(\phi_n)$ be a sequence of distinct characters on $A$ such that
the closure in $\Spec{A}$ of $\{\phi_1, \phi_2,\dots \}$
has no perfect subsets. Suppose that
$(\lambda_n) \in \ell_1$ is such that (with absolute operator norm convergence)
\[\sum_{n=1}^\infty{\lambda_n \phi_n} = 0\,.\]
Then all of the $\lambda_n$ must be zero.
\end{theorem}
\par
\Proof
The result follows immediately from Lemma 6 and Proposition 8.
\QED
\par
Note that, above, we have made heavy use of the fact the
$(\lambda_n) \in \ell_1$,
which enabled us to consider uniform algebras and annihilating measures. It
is not at all clear what happens if this condition on $(\lambda_n)$ is
removed, and some other form of convergence (such as conditional operator norm convergence or strong operator topology convergence)
is used when considering the series $\sum_{n=1}^\infty{\lambda_n \phi_n}$.
I do not know whether or not some or all of the results above still remain
valid in such cases.

I am grateful to Milne Anderson and John Wermer for useful discussions.


\begin{thebibliography}{6}

\bibitem{Bear}
Bear,~H.S., Complex function algebras, \emph{Trans. Amer. Math. Soc.} 90 (1959), 383--393.

\bibitem{BSZ}
Brown,~L., Shields,~A., Zeller,~K., On absolutely convergent exponential sums,
\emph{Trans. Amer. Math. Soc.} 96 (1960), 162--183.

\bibitem{Browder}
Browder,~A., \emph{Introduction to function algebras}, W. A. Benjamin, Inc., New
  York--Amsterdam, 1969.


\bibitem{dePaepe}
de Paepe, P. J.,
Essential function algebras with large Silov boundary,
\emph{Math. Scand.} 43 (1978), 325--328.

\bibitem{G}
Gamelin,~T.W., \emph{Uniform Algebras}, Prentice-Hall, Inc.,
Englewood Cliffs, N.~J., 1969.

\bibitem{R}
Rudin,~W., Continuous functions on compact spaces without perfect subsets.
\emph{Proc. Amer. Math. Soc.} 8 (1957), 39--42.

\bibitem{WO}
Wolff,~J., Sur les s\'eries $\sum\frac{A_k}{z-\alpha_k}$. \emph{C.~R.~Acad.
Sci. Paris} 173 (1921), 1057--1058, 1327--1328.
\end{thebibliography}
\end{document}